\documentclass[12pt]{amsart}
\usepackage{amsmath,amssymb,amsthm,mathtools}
\usepackage{mathrsfs,bm}
\usepackage{thmtools,thm-restate}
\usepackage{tikz-cd}
\usepackage{tikz}
\usepackage{enumitem}
\usepackage{graphicx}
\usetikzlibrary{matrix,arrows.meta,bending}
\usepackage{color}
\usepackage[numbers]{natbib}
\usepackage{booktabs}
\usepackage{longtable}
\usepackage{pgfplots}
\pgfplotsset{compat=1.14} 
\usetikzlibrary{arrows.meta,decorations.pathreplacing,positioning,calc} 

\usepackage[a4paper, left=3cm, right=3cm, top=3cm, bottom=3cm, footskip=15mm]{geometry}

\newtheorem{theorem}{Theorem}[section]

\newtheorem{lemma}[theorem]{Lemma}
\newtheorem{problem}[theorem]{Problem}

\newtheorem{corollary}[theorem]{Corollary}

\theoremstyle{definition}
\newtheorem{definition}[theorem]{Definition}

\theoremstyle{remark}

\numberwithin{equation}{section}

\makeindex

\usepackage{fancyhdr}
\usepackage{calc} 

\AtBeginDocument{%
  \setlength{\textwidth}{\dimexpr\paperwidth - 6cm\relax}%
  \setlength{\oddsidemargin}{\dimexpr 3cm - 1in\relax}%
  \setlength{\evensidemargin}{\dimexpr 3cm - 1in\relax}%
  \setlength{\marginparwidth}{2.5cm}%
}

\begin{document}

\title{On the distribution of addition chains}

\author{T. Agama}
\address{Department of Mathematics, African Institute for Mathematical science, Ghana
}
\email{theophilus@aims.edu.gh/emperordagama@yahoo.com}


\date{\today}



\begin{abstract}
Addition chains are a classical construction for fast exponentiation and related computation problems. In this paper, we study a chain for a fixed integer $n$ by decomposing each generator into a \emph{determiner} and a \emph{regulator} (gap). This viewpoint leads to explicit identities for two aggregate statistics of the chain: the sum of the determiners and the sum of the chain elements. We then derive the corresponding lower bounds by using the positivity of the regulators. In parallel, we establish an identity for the reciprocal sum of the chain, showing how the harmonic profile of the chain can also be written in terms of the same gap sequence. These identities provide a unified way to compare addition chains of the same target and length. The paper concludes with a balancing problem that asks for the chain(s) that minimize the difference between the arithmetic sum and the harmonic sum, together with a structural decomposition of that optimization objective.
\end{abstract}

\maketitle

\begingroup
  \setlength{\parskip}{6pt} 
  \tableofcontents
\endgroup

\section{Introduction}

Addition chains are one of the standard combinatorial devices for fast exponentiation: starting from $1$, each new term is obtained as the sum of two previous terms (repetition allowed). The length of the chain is the number of additions required to reach a fixed target and measures the number of multiplications required in the associated exponentiation procedure. The subject of an addition chain has a classical history, beginning with the early work of Alfred Brauer and later developed in the standard treatment by Knuth; modern surveys emphasize both the algorithmic relevance of the topic and the difficulty of producing the shortest chains in general \cite{brauer1939addition,knuth1981seminumerical,gordon1997survey,bergeron1994efficient}. In particular, while many constructive heuristics are known, the exact minimal length $\ell(n)$ is difficult to determine and is known only for relatively small values of $n$ \cite{gordon1997survey,bergeron1994efficient}.\\

The present paper does not focus on the shortest-chain problem itself. Instead, it studies the internal structure of a chosen addition chain by separating each generator into two parts: a determiner and a regulator (gap). This decomposition makes it possible to express aggregate quantities attached to the chain in terms of the regulator (gap) sequence alone. The first basic observation is a telescoping identity for the regulators (gaps), namely that their total sum is $n-1$. From this starting point, Theorem~\ref{identity} expresses the sum of determiners as an explicit weighted sum of the gaps, and Theorem~\ref{lengthidenity} turns that relation into a closed formula for the sum of chain elements. Each identity immediately yields a lower bound once the regulators are observed to be positive.\\

A second theme of the paper is harmonic rather than arithmetic. For the same chain, Theorem~\ref{harmonic length} analyzes the reciprocal sum $\sum 1/s_j$ and obtains an exact expansion that again depends only on the gap sequence. This gives a complementary description of the same chain data from a harmonic perspective and highlights how the shape of the gaps influences both large-scale and reciprocal statistics. The resulting identities are elementary in nature, but they make the dependence on the regulator (gap) sequence completely explicit.\\

The final section presents a balancing problem motivated by these two perspectives. Given a family of all addition chains that leads to a fixed target $n$, we ask for the chain(s) that minimize the difference between the arithmetic sum $\sum s_j$ and the harmonic sum $\sum 1/s_j$. The section also asks for a general decomposition of that objective into a term that depends only on the target $n$ and the length of the chain, together with correction terms determined by the gaps. This problem is intended as a natural optimization problem arising from the identities proved earlier in the paper.\\

\subsection{Organization of the paper} The paper is organized as follows. Section~2 introduces the terminology of addition chains, determiners, regulators (gaps), and truncated chains, and proves the basic telescoping lemma for regulators (gaps). Section~3 establishes the identities and lower bounds for the arithmetic (average) distribution of a chain. Section~4 derives the corresponding identities and lower bounds for the harmonic distribution. Section~5 formulates the balancing problem and collects the two main formulas proved in the preceding sections.

\section{Background}

\begin{definition}
Let $n\geq 3$. By an addition chain of length $k$ producing $n$, we mean the sequence of positive integers
\begin{align}
s_0=1,s_1=2,\ldots, s_{k}=n\nonumber
\end{align}
where each term $s_j=s_i+s_l$~($j\geq 1$) with repetition allowed, and with the corresponding sequence of partitions
\begin{align}
2=1+1,\ldots,s_{k-1}=a_{k-1}+r_{k-1},s_{k}=a_{k}+r_{k}=n\nonumber
\end{align}
where $a_{i+1}=a_{i}+r_{i}$ and $a_{i+1}=s_i$ for $1\leq i\leq k$. We call the partition $a_{i}+r_{i}$ the $i^{th}$ \emph{generator} of the chain for $2\leq i\leq k$. We call $a_{i}$ the \emph{determiner} and $r_{i}$ the \emph{regulator} of the $i^{th}$ generator of the chain. We call the sequence $(r_{i})$ the regulators (gap sequence) of the addition chain and $(a_i)$ the determiners of the chain for $2\leq i \leq k$. We call the subsequence $(s_{j_{m}})$ for $2\leq j\leq k$ and $1\leq m\leq t\leq k$ a truncated addition chain producing $n$.
\end{definition}

\begin{lemma}\label{first}
Let $s_0=1,s_1=2,\ldots,s_k=n$ be an addition chain producing $n\geq 3$ with associated generators
\begin{align}
2=1+1,\ldots,s_{k-1}=a_{k-1}+r_{k-1},s_{k}=a_{k}+r_{k}=n.\nonumber
\end{align} 
We have the following relation for the regulators (gaps) 
\begin{align}
\sum \limits_{j=1}^{k}r_j=n-1.\nonumber
\end{align}
\end{lemma}

\begin{proof}
We can write $r_{k}=n-a_{k}$. It follows that 
\begin{align}
r_{k}+r_{k-1}&=n-a_{k}+r_{k-1}\nonumber \\&=n-(a_{k-1}+r_{k-1})+r_{k-1}\nonumber \\&=n-a_{k-1}.\nonumber
\end{align}
We deduce
\begin{align}
r_{k}+r_{k-1}+r_{k-2}&=n-a_{k-1}+r_{k-2}\nonumber \\&=n-(a_{k-2}+r_{k-2})+r_{k-2}\nonumber \\&=n-a_{k-2}.\nonumber
\end{align}
Iterating downward in this manner, the relation follows.
\end{proof}

\section{The average distribution}

Let $s_0=1,\ldots,s_h=n$ be an addition chain that leads to $n$ and of length $h$. We define $s_j-s_{j-1}:=r_{j}$ as the absolute difference between the consecutive terms in the addition chain. The goal of this section is to prove an identity of the form 
$$
\sum \limits_{j=0}^{h}s_j=L(n,h)+\sum \limits_{j=1}^{h-1}R(n,h,r_j).
$$

\begin{theorem}\label{identity}
Let $s_0=1,s_1=2,\ldots,s_{\delta(n)}$ be an addition chain producing $n\geq 3$ of length $\delta(n)$ with associated generators 
\begin{align}
s_1=a_1+r_1,s_2=a_2+r_2,\ldots,s_{\delta(n)}=a_{\delta(n)}+r_{\delta(n)}=n.\nonumber
\end{align}
We have
\begin{align}
\sum \limits_{j=1}^{\delta(n)}a_j=\delta(n)+\sum \limits_{j=1}^{\delta(n)-1}(\delta(n)-j)r_j.\nonumber
\end{align}
\end{theorem}

\begin{proof}
This identity is easily established by induction on an iteration on the sum. We first observe
\begin{align}
a_1+a_2&=a_1+a_1+r_1\nonumber \\&=2a_1+r_1\nonumber \\&=2+r_1.\nonumber \quad (a_1=s_0=1)
\end{align}
Again, it follows that 
\begin{align}
a_1+a_2+a_3&=2a_1+r_1+a_3\nonumber \\&=3a_1+2r_1+r_2\nonumber \\&=3+2r_1+r_2.\nonumber \quad (a_1=s_0:=1)
\end{align}
By induction, we can write
\begin{align}
\sum \limits_{j=1}^{\delta(n)}a_j&=\delta(n)+(\delta(n)-1)r_1+(\delta(n)-2)r_2+\cdots+r_{\delta(n)-1}\nonumber
\end{align}
 establishing the identity.
\end{proof}
\bigskip

The above identity suggests that the partial sums of the determiners $a_j$ of the generators of any addition chain can be written as a linear combination of the regulators (gaps) $r_j$.

\begin{theorem}\label{lengthidenity}
Let $s_0=1,s_1=2,\ldots,s_{\delta(n)}=n$ be an addition chain producing $n\geq 3$ of length $\delta(n)$ with associated generators 
\begin{align}
s_1=2=1+1,s_2=a_2+r_2,\ldots,s_{\delta(n)}=a_{\delta(n)}+r_{\delta(n)}=n.\nonumber
\end{align}
We have
\begin{align}
\sum \limits_{j=1}^{\delta(n)}s_j=(n-1)+\delta(n)+\sum \limits_{j=1}^{\delta(n)-1}(\delta(n)-j)r_j.\nonumber
\end{align}
\end{theorem}

\begin{proof}
By Theorem \ref{identity}, we write
\begin{align}
    \sum \limits_{j=1}^{\delta(n)}s_j&=\sum \limits_{j=1}^{\delta(n)}a_j+\sum \limits_{j=1}^{\delta(n)}r_j\nonumber \\&=(n-1)+\delta(n)+\sum \limits_{j=1}^{\delta(n)-1}(\delta(n)-j)r_j.\nonumber
\end{align}
\end{proof}
\bigskip

Here, we deduce a lower bound for the partial sum $\sum s_j$ for an addition chain $s_0=1,s_1=2,\ldots,s_{h}=n$ that leads to $n$.

\begin{corollary}
    Let $s_0=1,s_1=2,\ldots,s_{\delta(n)}=n$ be an addition chain that leads to $n$ and of length $\delta(n)$. We have 
    $$
    \sum \limits_{j=1}^{\delta(n)}s_j\geq (n-1)+\frac{\delta(n)(\delta(n)+1)}{2}.
    $$
\end{corollary}

\begin{proof}
    We observe $r_j\geq 1$ and deduce (by Theorem \ref{lengthidenity})
    \begin{align}
        \sum \limits_{j=1}^{\delta(n)}s_j&=(n-1)+\delta(n)+\sum \limits_{j=1}^{\delta(n)-1}(\delta(n)-j)r_j\nonumber \\&\geq (n-1)+\delta(n)+\sum \limits_{j=1}^{\delta(n)-1}(\delta(n)-j)\nonumber \\&=(n-1)+\delta(n)+\frac{\delta(n)(\delta(n)-1)}{2}.\nonumber
    \end{align}
\end{proof}

\section{The Harmonic distribution}

Let $s_0=1,\ldots,s_h=n$ be an addition chain that leads to $n$ and of length $h$. We define $s_j-s_{j-1}:=r_{j}$ as the absolute difference between the consecutive terms in the addition chain. The goal of this section is to prove an identity of the form 
$$
\sum \limits_{j=0}^{h}\frac{1}{s_j}=L(n,h)+\sum \limits_{j=1}^{h-1}R(n,h,r_j).
$$

\begin{theorem}\label{harmonic length}
Let $s_0=1,s_1=2,\ldots,s_{\delta(n)}=n$ be an addition chain producing $n$ and of length $\delta(n)$, with an associated sequence of generators 
\begin{align}
s_1=1+1,s_2=a_2+r_2,\ldots,s_{\delta(n)}=a_{\delta(n)}+r_{\delta(n)}=n.\nonumber
\end{align}
We have
$$
\sum \limits_{l=0}^{\delta(n)}\frac{1}{s_l}=\frac{3}{2}+\frac{\delta(n)-1}{n}+\sum \limits_{l=2}^{\delta(n)-1}\sum \limits_{v=1}^{\infty}\frac{1}{n^{v+1}}\bigg(\sum \limits_{j=l+1}^{\delta(n)}r_j\bigg)^{v}
$$
where $\sum \limits_{j=l}^{\delta(n)}r_j<n-1$ is the sum of the regulators (gaps) in the chain generator for each $2\leq l\leq \delta(n)$.
\end{theorem}

\begin{proof}
We consider an addition chain 
$$
s_0=1,s_1=2,\ldots,s_{\delta(n)}=n
$$ 
producing $n$ and of length $\delta(n)$, with the associated sequence of generators 
\begin{align}
s_1=1+1,s_2=a_2+r_2,\ldots,s_{\delta(n)}=a_{\delta(n)}+r_{\delta(n)}=n\nonumber
\end{align} 
and let $(a_j)$ and $(r_j)$ be the sequence of determiners and regulators (gaps), respectively, in the chain. We make the following observations: $s_{\delta(n)-1}=a_{\delta(n)}=a_{\delta(n)-1}+r_{\delta(n)-1}=s_{\delta(n)-2}+r_{\delta(n)-1}=a_{\delta(n)-2}+r_{\delta(n)-2}+r_{\delta(n)-1}=\cdots=1+\sum \limits_{j=1}^{\delta(n)-1}r_j=n-r_{\delta(n)}$, where we have used Lemma \ref{first}. Similarly, we can write $a_{\delta(n)-1}=1+\sum \limits_{j=1}^{\delta(n)-2}=n-r_{\delta(n)}-r_{\delta(n)-1}$. By induction, we can write 
$$
a_l=n-\sum \limits_{j=l}^{\delta(n)}r_j
$$ 
for each $3\leq l\leq \delta(n)$. We observe that 
$$
\sum \limits_{l=0}^{\delta(n)}\frac{1}{s_l}=\frac{3}{2}+\sum \limits_{l=3}^{\delta(n)}\frac{1}{a_l}+\frac{1}{n}=\frac{3}{2}+\sum \limits_{l=2}^{\delta(n)-1}\frac{1}{s_l}+\frac{1}{n}.
$$ 
We now analyze the latter sum of the right-hand side that involves the \emph{determiners} of the addition chain. We can write 
$$
\sum \limits_{l=2}^{\delta(n)-1}\frac{1}{s_l}=\sum \limits_{l=2}^{\delta(n)-1}\frac{1}{n-\sum \limits_{i=l+1}^{\delta(n)}r_i}
$$ 
which can be recast as
$$
\sum \limits_{l=2}^{\delta(n)-1}\frac{1}{s_l}=\sum \limits_{l=2}^{\delta(n)-1}\frac{1}{n}+\sum \limits_{l=2}^{\delta(n)-1}\sum \limits_{v=1}^{\infty}\frac{1}{n^{v+1}}\bigg(\sum \limits_{i=l+1}^{\delta(n)}r_i\bigg)^v
$$ 
with $\sum \limits_{i=l+1}^{\delta(n)}r_i<n-1$ for each $2\leq l\leq \delta(n)$ by Lemma \ref{first}. It follows that 
$$
\sum \limits_{l=0}^{\delta(n)}\frac{1}{s_{l}}=\frac{3}{2}+\frac{\delta(n)-2}{n}+\sum \limits_{l=2}^{\delta(n)-1}\sum \limits_{v=1}^{\infty}\frac{1}{n^{v+1}}\bigg(\sum \limits_{i=l+1}^{\delta(n)}r_i\bigg)^v+\frac{1}{n}
$$
where $\sum \limits_{i=l}^{\delta(n)}r_i<n-1$ by Lemma \ref{first} for each $2\leq l\leq \delta(n)$. This completes the proof of the claimed identity.
\end{proof}
\bigskip

Here, we deduce a lower bound for the partial sum $\sum \frac{1}{s_j}$ for an addition chain $s_0=1,s_1=2,\ldots,s_{h}=n$ that leads to $n$.

\begin{corollary}
    Let $s_0=1,s_1=2,\ldots,s_{\delta(n)}=n$ be an addition chain that leads to $n$ and of length $\delta(n)$. We have
    $$
    \sum \limits_{l=0}^{\delta(n)}\frac{1}{s_l}\geq \frac{3}{2}+\frac{\delta(n)-1}{n}+\sum \limits_{l=2}^{\delta(n)-1}\sum \limits_{v=1}^{\infty}\frac{1}{n^{v+1}}\bigg(\delta(n)-l\bigg)^{v}.
    $$
\end{corollary}

\begin{proof}
    We observe $r_j\geq 1$ and deduce (by Theorem \ref{harmonic length})
    \begin{align}
        \sum \limits_{l=0}^{\delta(n)}\frac{1}{s_l}&=\frac{3}{2}+\frac{\delta(n)-1}{n}+\sum \limits_{l=2}^{\delta(n)-1}\sum \limits_{v=1}^{\infty}\frac{1}{n^{v+1}}\bigg(\sum \limits_{j=l+1}^{\delta(n)}r_j\bigg)^{v}\nonumber \\&\geq \frac{3}{2}+\frac{\delta(n)-1}{n}+\sum \limits_{l=2}^{\delta(n)-1}\sum \limits_{v=1}^{\infty}\frac{1}{n^{v+1}}\bigg(\sum \limits_{j=l+1}^{\delta(n)}1\bigg)^{v}\nonumber \\&=\frac{3}{2}+\frac{\delta(n)-1}{n}+\sum \limits_{l=2}^{\delta(n)-1}\sum \limits_{v=1}^{\infty}\frac{1}{n^{v+1}}\bigg(\delta(n)-l\bigg)^{v}.\nonumber
    \end{align}
\end{proof}

\section{A balancing problem}

In this section, we pose a natural optimization problem that naturally emerges from the study of the distribution of an addition chain that leads to a fixed target. The study has been twofold: the \emph{average} and \emph{harmonic} distribution of an addition chain that leads to a fixed target and of a given length. Precisely, given an addition chain 
$$
s_0=1,s_1=2,\ldots,s_{\delta(n)}=n
$$
with $r_i=s_i-s_{i-1}$ for $1\leq i\leq h$, we have proved the following identities: 

\begin{align}
\sum \limits_{j=1}^{\delta(n)}s_j=(n-1)+\delta(n)+\sum \limits_{j=1}^{\delta(n)-1}(\delta(n)-j)r_j\nonumber
\end{align}
and 
$$
\sum \limits_{l=0}^{\delta(n)}\frac{1}{s_l}=\frac{3}{2}+\frac{\delta(n)-1}{n}+\sum \limits_{l=2}^{\delta(n)-1}\sum \limits_{v=1}^{\infty}\frac{1}{n^{v+1}}\bigg(\sum \limits_{j=l+1}^{\delta(n)}r_j\bigg)^{v}.
$$
\bigskip

However, the question of how much the reciprocal sum $\sum \frac{1}{s}$ deviates from the partial sum $\sum s$ may be an interesting problem for investigation. In other words, it is an interesting problem to determine from among the family of all chains that leads to a fixed target the chain that minimizes the difference
$$
\sum s_i-\sum \frac{1}{s_i}.
$$
In particular, we pose the following problem:

\begin{problem}
    Let $n$ be a fixed positive integer and $\mathcal{L}_n$ denote the family of all addition chains that leads to $n$. Determine the chain(s) in $\mathcal{L}_n$ that minimizes the difference 
    $$
    \sum s_i-\sum \frac{1}{s_i}
    $$
    where $s_i$ are the terms in the chosen addition chain. Furthermore, deduce an expression of the form 
    $$
    \sum s_i-\sum \frac{1}{s_i}=G(n,\mathrm{length})+\sum \limits_{i=1}^{\mathrm{length}-1}H_i(n,\mathrm{length, gap}).
    $$
\end{problem}
\footnote{
\par
.}%

\bibliographystyle{amsplain}

\end{document}